\documentclass[reqno, 11pt]{amsart}

\newtheorem{theorem}{Theorem}[section]

\theoremstyle{definition}

\theoremstyle{remark}
\newtheorem{remark}[theorem]{Remark}

\numberwithin{equation}{section}

\usepackage{bbm}
\usepackage{euscript}
\usepackage{pb-diagram}
\usepackage{lamsarrow}
\usepackage{pb-lams}
\usepackage{amsmath}
\usepackage{amsthm}
\usepackage{epsfig}
\usepackage{amssymb}
\usepackage{pifont}
\usepackage{amsbsy}

\newskip\aline \newskip\halfaline
\def\skipaline{\vskip\aline}

\def\qedbox{$\rlap{$\sqcap$}\sqcup$}
\def\qed{\nobreak\hfill\penalty250 \hbox{}\nobreak\hfill\qedbox\skipaline}

\newcommand\bR{{\mathbb R}}

\newcommand\bZ{{\mathbb Z}}

\newcommand{\m}{\EuScript{M}}

\newcommand{\p}{\EuScript{P}}

\newcommand{\s}{\EuScript{S}}

\newcommand{\ra}{\rightarrow}

\def\inpr{\mathbin{\hbox to 6pt{\vrule height0.4pt width5pt depth0pt \kern-.4pt \vrule height6pt width0.4pt depth0pt\hss}}}

\newcommand{\vfi}{{\varphi}}

\newcommand{\pa}{\partial}

\begin{document}

\title{Morse functions statistics}
\date{April 20, 2006.}

\author{Liviu I. Nicolaescu}
\address{Department of Mathematics, University of Notre Dame, Notre Dame, IN 46556-4618.}
\email{nicolaescu.1@nd.edu} \urladdr{http://www.nd.edu/~lnicolae/}
\thanks{This work was partially supported by NSF grant DMS-0303601.}

\begin{abstract} We answer a question of V.I. Arnold concerning the growth rate of the number of Morse functions on the two sphere.
\end{abstract}

\maketitle

\section{Introduction}
\setcounter{equation}{0}

We  are interested in \emph{excellent} Morse functions $f: S^2\ra \bR$, where the attribute  excellent signifies that no two critical points lie on the same level set of $f$. Two such Morse functions $f_0,f_1$ are called geometrically  equivalent  if there exist orientation preserving diffeomorphisms  $R: S^2\ra S^2$ and $L: \bR\ra \bR$ such that   $f_1=  L\circ f_0\circ R^{-1}$.  We denote  by $ g(n)$ the number of equivalence classes of Morse functions  with $2n+2$ critical points.  Arnold    suggested in \cite{Ar} that
\begin{equation}
\lim_{n\ra \infty}\frac{\log g(n)}{n\log n}=2.
\label{eq: asympt}
\end{equation}
The goal  of this note is to  establish the validity of Arnold's prediction.

\smallskip

\noindent {\bf Acknowledgment.}   I want to thank Francesca Aicardi  for  drawing my attention to Arnold's question.

\section{Some   background on the number of Morse functions}
\setcounter{equation}{0}

We define
\[
h(n):= \frac{g(n)}{(2n+1)!},\;\;\xi(\theta):=\sum_{n\geq 0}  h(n)\theta^{2n+1}.
\]
In \cite{N} we have embedded $h(n)$ in a $2$-parameter family
\[
(x,y)\longmapsto  \hat{H}(x,y),\;\; x,y\in \bZ_{\geq 0},\;\;  h(n) =\hat{H}(0,n)
\]
which satisfies a nonlinear recurrence relation, \cite[\S 8]{N}.

\smallskip

\noindent {\bf A.} $x>0$.
\[
(x+2y+1)\hat{H}(x,y) - (x+1) \hat{H}(x+1,y-1)
\]
\[
=\frac{x+1}{2} \hat{H}(x-1,y)+\frac{x+1}{2}\sum_{(x_1,y_1)\in R_{x,y-1}} \hat{H}(x_1,y_1)\hat{H}(\bar{x}_1,\bar{y}_1),
\]
where
\[
R_{x,y-1}=\{ (a,b)\in \bZ^2;\;\; 0\leq a \leq x,\;\; 0\leq b\leq y-1\},
\]
and for every $(a,b)\in R_{x,y-1}$ we denoted by $(\bar{a},\bar{b})$ the  symmetric  of $(a,b)$  with respect to the center of  the rectangle $R_{x,y-1}$.

\smallskip

\noindent {\bf B.} $x=0$.
\[
(2y+1)\hat{H}(0, y)-\hat{H}(1,y-1)=\frac{1}{2}\sum_{y_1=0}^{y-1} \hat{H}(0,y_1)\hat{H}(0,y-1-y_1).
\]
Observe that if  we let $y=0$ in {\bf A} we deduce
\[
\hat{H}(x,0)=\frac{1}{2} \hat{H}(x-1,0)
\]
so that $\hat{H}(x,0)=2^{-x}$.

In \cite{N} we proved that these recurrence relations imply that the  function
\[
\xi(u,v)=\sum_{x,y\geq 0}\hat{H}(x,y) u^xv^{x+2y+1}
\]
 satisfies the   quasilinear pde
 \[
 -\bigl(1+u\xi+\frac{u^2}{2}\,\bigr)\pa_u\xi+\pa_v\xi= \bigl(\frac{1}{2}\xi^2+u\xi+1),\;\;\xi(u,0)=0,
 \]
 and  inverse function $\xi(0,\theta) =\xi \longmapsto \theta$ is defined by the elliptic integral
\begin{equation}
\theta=\int_0^\xi\frac{dt}{\sqrt{t^4/4-t^2+2\xi t+1}} .
\label{eq: ell}
\end{equation}

\section{Proof of the asymptotic estimate}

Using the recurrence formula {\bf B}    we deduce that for every $n\geq 1$ we have
\[
(2n+1) h(n) \geq \frac{1}{2}\sum_{k=0}^{n-1} h(k)h(n-1-k)
\]
We multiply this  equality with $t^{2n}$ and we deduce
\[
\sum_{n\geq 1} (2n+1) h(n) t^{2n} \geq  \frac{1}{2}\sum_{n\geq 1}\Biggl( \sum_{k=0}^{n-1} h(k)h(n-1-k)\Biggr) t^{2n}
\]
$(g(0)=1)$
\[
\Longleftrightarrow \frac{d\xi}{dt} \geq 1+\frac{1}{2} \xi^2.
\]
This implies that the Taylor coefficients of $\xi$ are bounded from below by the Taylor coefficients of the solution of the initial value problem
\[
\frac{du}{dt}=1+\frac{1}{2}u^2,\;\;u(0)=\xi(0)=0.
\]
 The  last  initial value problem  can be solved by separation of variables
\[
\frac{du}{1+\frac{u^2}{2}} =dt\Longrightarrow u=\sqrt{2} \tan (t/\sqrt{2}).
 \]
The  function $\tan$ has the Taylor  series (see \cite[\S 1.41]{GR})
\[
\tan x =\sum_{k=1}^\infty \frac{2^{2k}(2^{2k}-1)|B_{2k}|}{(2k)!} x^{2k-1},
\]
where  $B_{n}$ denote the Bernoulli numbers generated by
\[
\frac{t}{e^t-1}=\sum_{n=0}^\infty B_n \frac{t^n}{n!}.
\]
The  Bernoulli numbers have the asymptotic behavior, \cite[Sec. 6.2]{B}
\[
|B_{2k}| \sim  \frac{2(2k)! }{(4\pi^2)^k}.
\]
If $T_k$ denotes the coefficient  of $x^{2k+1}$ in  $\tan(x)$ we deduce that
\[
T_k= \frac{2^{2k+2}(2^{2k+2}-1)|B_{2k+2}|}{(2k+2)!}\sim \frac{2^{2k+3}(2^{2k+2}-1)}{(4\pi^2)^{k+1}}.
\]
Thus the coefficient $u_k$ of $t^{2k+1}$ in $\sqrt{2}\tan(t/\sqrt{2})$ has the asymptotic behavior
\[
u_k\sim \frac{1}{2^k} \frac{2^{2k+3}(2^{2k+2}-1)}{(4\pi^2)^{k+1}}= \frac{2^{k+3}(2^{2k+2}-1)}{(4\pi^2)^{k+1}}.
\]
We deduce that
\begin{equation*}
g(k) > (2k+1)! \frac{2^{k+3}(2^{2k+2}-1)}{(4\pi^2)^{k+1}}(1+o(1)\,),\;\;\mbox{as $k\ra \infty$}.
\tag{$\dag$}
\label{tag: dag}
\end{equation*}
Let us produce  upper bounds for $g(n)$.  We will give a  combinatorial argument showing that
\[
g(n)\leq (2n+1)!  C_n,
\]
where $C_n=\frac{1}{n+1}\binom{2n}{n}$ is the $n$-th Catalan number.

\begin{figure}[ht]
\centerline{\epsfig{figure=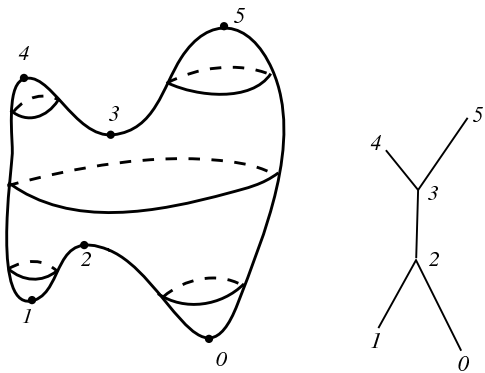,height=1.9in,width=2.5in}}
\caption{\sl Associating a tree to a Morse function on $S^2$.}
\label{fig: 1}
\end{figure}

As  explained in \cite{Ar,N}, a  geometric equivalence class  of a Morse function  on $S^2$ with $2n+2$ critical points is completely described by a certain labelled tree, dubbed  Morse tree in \cite{N} (see Figure \ref{fig: 1}, where the Morse function is the height function).  For the reader's convenience we recall that a Morse tree (with $2n+2$ vertices)  is tree with vertices  labelled   with labels  $\{0,1,\cdots\}$ with the following two properties.

\smallskip

\noindent $\bullet$ Any vertex  has either one neighbor, or exactly three neighbors, in which case the vertex is called a node.

\noindent $\bullet$ Every node has at least one neighbor with a    higher label, and at least one neighbor with a lower label.

\smallskip

We will produce an injection from the set $\m_n$ of Morse functions with $2n+2$ critical points to the set $\p_n\times S_{2n+1}$ where $\p_n$ denotes the set of  Planted, Trivalent, Planar  Trees (PTPT) with  $2n+2$-vertices, and $S_{2n+1}$ denotes the group of permutations  of $2n+1$-objects.

As explained in \cite[Prop. 6.1]{N}, to a Morse tree we can  canonically associate a   PTPT with $2n+2$-vertices. The number of such PTPT-s is $C_n$,  \cite[Exer. 6.19.f, p.220]{St2}.  The tree in Figure \ref{fig: 1} is already a $PTPT$.

\begin{figure}[ht]
\centerline{\epsfig{figure=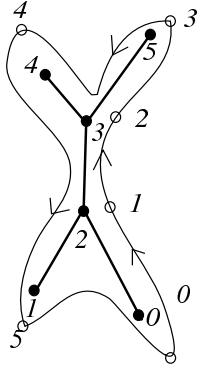,height=1.9in,width=1in}}
\caption{\sl Labelling the vertices of a $PTPT$.}
\label{fig: 2}
\end{figure}

The   non-root vertices  of such a tree can be  labelled  in a canonical way  with labels $\{1,2,\cdots, 2n+1\}$ (see the explanation in  \cite[Figure 5.14, p. 34]{St2}).  More precisely, consider a  very thin tubular neighborhood $N$ of such a tree  in the plane.  Its boundary   is a circle.    To label the vertices,   walk along $\pa N$   in a counterclockwise fashion  and  label the non-root vertices   in the order they were first encountered (such a walk  passes three times near each node).  In Figure \ref{fig: 2},  this labelling is indicated  along the points marked $\circ$. The Morse function then   defines   another bijection from the set of non-root vertices  to   the same label set. In Figure \ref{fig: 2} this labelling is indicated  along the vertices marked $\bullet$.

 We have thus associated  to a Morse tree   a pair, $(T,\vfi)$, where $T$ is a PTPT and $\vfi$ is a permutation of its non-root vertices.  In Figure \ref{fig: 2} this permutation is
 \[
 1\ra 2,\;\;2\ra 3,\;\;3\ra 5,\;\;4\ra 4,\;\;5\ra 1.
 \]
  The  Morse tree is uniquely determined by this pair. We deduce that
\[
g(n)=\# \m_n \leq \# \p_n \times \# S_{2n+1}= C_n (2n+1)!
\]
\[
=\frac{(2n)!}{(n+1)!n!}(2n+1)! =\frac{2\cdot 4\cdots (2n) \cdot 1\cdot 3\cdot 5\cdots (2n-1)}{n!\cdot n!}\cdot \frac{(2n+1)!}{n+1}.
\]
Hence
\begin{equation*}
g(n) <2^n \cdot \frac{(2n+1)!}{n+1}\prod_{k=0}^{n-1}\frac{2k-1}{k+1}\leq  \frac{2^{2n}}{(n+1)} (2n+1)!.
\tag{$\ddag$}
\label{tag: ddag}
\end{equation*}
The   estimates  (\ref{tag: dag}) and (\ref{tag: ddag}) coupled with Stirling's formula  show that
\[
\lim_{n\ra \infty}\frac{\log g(n)}{n\log n}=2,
\]
which is Arnold's prediction, (\ref{eq: asympt}).  \qed

\begin{remark}  (a) Numerical  experiments suggest that
\[
g(n)  < (2n+1)!.
\]
Is it possible to give a purely combinatorial proof of this inequality?

\noindent (b) It would be interesting to have a more refined  asymptotic  estimate for   $g(n)$ of the form
\[
\log g(n)= 2n\log n + r_n,\;\; r_n=  an + b \log n + c + O(n^{-1}),\;\;a,b,c\in \bR.
\]
Here are the results of some numerical  experiments.

The refined Stirling's formula
\[
\log(2n+1)!= (2n+3/2)\log(2n+1)-2n-1+\frac{1}{2}\log(2\pi)+O(n^{-1})
\]
implies that
\[
\log h(n) = \log g(n)-\log(2n+1)!
\]
\[
= 2n\log n + r_n -  (2n+3/2)\log(2n+1)+2n+1-\frac{1}{2}\log(2\pi)+O(n^{-1})
\]
\[
=r_n +2n\Bigl( 1+\log\frac{n}{2n+1}\Bigr)-3/2\log(2n+1)+1-\frac{1}{2}\log(2\pi)+O(n^{-1}).
\]
Hence
\[
r_n =\underbrace{\log h(n)- 2n\Bigl( 1+\log\frac{n}{2n+1}\Bigr) +3/2\log(2n+1)-1 +\frac{1}{2}\log(2\pi)}_{\delta_n}+O(n^{-1}).
\]
We deduce that
\[
\frac{r_n}{n}=\frac{\delta_n}{n}+O(n^{-2}).
\]
Here are the  results of some numerical  experiments.
\[
\begin{tabular}{||r|r||}\hline
$n$  & $\delta_n/n$\\ \hline
$10$ & $-0.634$\\
$20$ & $-0.750$\\
$30$ & $-0.790$ \\
$40$ & $-0.811$\\
$50$ & $-0.824$\\
$100$& $-0.849$\\
$150$ & $-0.858$\\
$200$ & $-0.862$ \\ \hline\hline
\end{tabular}
\]
This suggests $ a\approx -0.8\cdots$.
\end{remark}


\begin{thebibliography}{XXXXXX}

\bibitem{Ar} V.I. Arnold: {\sl Smooth functions statistics}, preprint,

\verb+ http://www.institut.math.jussieu.fr/seminaires/singularites/Arnold.html +

\bibitem{B} E.A. Bender: {\sl Asymptotic enumeration}, SIAM Rev., {\bf 16}(1974), 485-515.

%\bibitem{MM}  A. Meir, J.W. Moon: {\sl On an asymptotic method in enumeration}, J. Comb. Theory, Series A, {\bf 51}(1989), 77-89.

\bibitem{GR} I.S. Gradshteyn,  I.M. Ryzhik, {\sl Table of Integrals, Series and Products}, Academic Press, 2000.

\bibitem{N}  L.I. Nicolaescu: {\sl Counting Morse functions on the $2$-sphere}, preprint, {\sf math.GT/0512496}.

%\bibitem{O} A.M. Odlyzko: {\sl  Asymptotic enumeration methods},    {\sl Handbook of Combinatorics, vol 2}, p.1063-1230,  Elsevier-MIT Press, 1995.

\bibitem{St2}  R.P. Stanley: {\sl Enumerative Combinatorics. Volume II}, Cambridge Stud. Adv. Math., vol. 62, Cambridge University Press, 1999.




\end{thebibliography}
\end{document}